\documentclass{amsart}
\usepackage{amsmath, amssymb, euscript,  enumerate}

\newtheorem{thm}{Theorem}[section]
\newtheorem{cor}[thm]{Corollary}
\newtheorem{lem}[thm]{Lemma}
\newtheorem{claim}[thm]{Claim}
\newtheorem{remark}[thm]{Remark}
\newtheorem{prop}[thm]{Proposition}

\theoremstyle{remark}
\newtheorem{rem}[thm]{Remark}
\numberwithin{equation}{section}
\newcommand{\R}{\mathbb R}

\newcommand{\tu}{\tilde u}

\newcommand{\dv}{\mathrm{div}}
\newcommand{\sign}{\mathrm{sign}}

\def\XXint#1#2#3{{\setbox0=\hbox{$#1{#2#3}{\int}$}
     \vcenter{\hbox{$#2#3$}}\kern-.5\wd0}}

\begin{document}
\title{On the extinction profile  of solutions to
fast-diffusion}

\author{Panagiota Daskalopoulos$^*$}

\address{Department of
Mathematics, Columbia University, New York,
 USA}
\email{pdaskalo@math.columbia.edu}

\thanks{$*:$ Partially supported
by NSF grant 0102252}

\author{Natasa Sesum}
\address{Department of Mathematics, Columbia University, New York,
USA}
\email{natasas@math.columbia.edu}

\begin{abstract}
We study the extinction behavior   of solutions to the fast diffusion equation $u_t = \Delta u^m$ on  $\R^N\times (0,T)$, in the
range of exponents $m \in (0, \frac{N-2}{N})$, $N > 2$.  We show
that if the initial data $u_0$ is trapped in between two Barenblatt solutions vanishing at time $T$,
then the vanishing behaviour of $u$ at $T$ is given by a Barenblatt solution.
We also give an example showing that for such a behavior  the bound from above  by a Barenblatt solution $B$
(vanishing at $T$)  is crucial:  we construct a class of solutions $u$ with initial data
$u_0 = B\, (1 + o(1))$, near  $ |x| >> 1$,  which live  longer than $B$
and change   behaviour at $T$. The  behavior of such solutions is governed by $B(\cdot,t)$ up to $T$,  while for $t >T$
the solutions   become  integrable and exhibit  a different vanishing profile. For
the Yamabe flow ($m = \frac{N-2}{N+2}$) the above means that these  solutions $u$  develop a singularity at time $T$, when the Barenblatt solution disappears, and at $t >T$  they
immediately smoothen  up and exhibit   the vanishing profile of a sphere.

In the appendix we show how we remove the assumption on the bound
on $u_0$ by a Barenblatt from below.
\end{abstract}

\maketitle

\begin{section}{Introduction}\label{sec-intro}

We consider the Cauchy problem for the fast diffusion equation
\begin{equation}
\begin{cases} \label{eqn-u}
u_t  = \Delta u^m  & \mbox{in} \,\,  \R^N
\times (0,T)\\
u(x,0) = u_0(x) & x \in \R^N
\end{cases}
\end{equation}
in the range of exponents $0 < m < (N-2)/N$, in dimensions $N \geq 3$. The initial data $u_0$ is assumed to be non-negative and locally integrable.

Equation \eqref{eqn-u} arises as a model of  various diffusion processes.
It is found  in plasma physics and in particular  as the Okuda-Dawson low,  when   $m=1/2$. Also, King studies \eqref{eqn-u} in a model of diffusion of impurities in silicon \cite{King2}.

When $m=(N-2)/(N+2)$ equation \eqref{eqn-u} is equivalent  to
the evolution of the   conformally flat  metric $g_{ij} = u^{\frac{4}{N+2}}\, dx_idx_j$ by
the Yamabe flow
$$ \frac{\partial g}{\partial t}  = - R \,  g $$
where $R$ denotes the scalar curvature with respect to  the metric $g$. The equivalence easily follows  from the observation that the
conformal metric $g_{ij} = u^{\frac{4}{N+2}}\, dx_idx_j$ has scalar curvature
$$R = - \frac{\Delta  u^{\frac{N-2}{N+2}}}{u}.$$
For an introduction to the Yamabe flow see in \cite{Ye}.

Our goal in this paper is to study the vanishing behavior of solutions
to the Cauchy problem \eqref{eqn-u},  under the assumption that
the initial data $u_0$ satisfies the growth condition
\begin{equation}\label{eqn-igc1}
u_0(x) = \left ( \frac{C}{|x|^2}  \right )^{\frac 1{1-m}} (1 + o(1)), \qquad \mbox{as} \,\, |x| \to \infty.
\end{equation}
The self-similar Barenblatt solutions of equation \eqref{eqn-u} given explicitly  by
\begin{equation}\label{eqn-baren}
B_k(x,t) = \left (\frac{C^*\, (T-t)}{k \, (T-t)^{2\gamma}+ |x|^2} \right )^{\frac{1}{1-m}}
\end{equation}
with
\begin{equation}\label{eqn-parameters}
\beta = \frac{N}{N-2-Nm} \qquad \mbox {and} \qquad \gamma = -\frac{\beta}{N}
\end{equation}
satisfy the growth condition \eqref{eqn-igc1}. The constant
$C^*$ depends only on $m$ and $N$ and is given explicitly by
\begin{equation}\label{eqn-C*}
C^*= \frac{2m\,(N-2-m\, N)}{1-m} .
\end{equation}

We will assume in the first part of this paper that the  initial condition  $u_0$ is trapped in between two
Barenblatt  solutions, i.e.,
\begin{equation}
\label{eqn-trapped}
\left (\frac{C^*\, T}{k_1 + |x|^2}\right )^{\frac{1}{1-m}} \le u_0(x) \le \left (\frac{C^*\, T}{k_2 + |x|^2} \right )^{\frac{1}{1-m}}
\end{equation}
for some constants $k_1 > k_2 >0$. As a direct consequence of
the maximum principle we then have
\begin{equation}
\label{eqn-trapped2}
B_{k_1}(x,t) \leq u(x,t) \leq B_{k_2}(x,t) \qquad \mbox{for} \,\, 0 < t < T.
\end{equation}
In particular, $u$ vanishes at time $T$.
We will show  in the first part of this paper  that the vanishing profile of
$u$ is given by a Barenblatt solution.

%Assume that $u$ is a solution of \eqref{eqn-u} that satisfies condition
%\eqref{eqn-trapped2}. We will show  that the vanishing profile of
%$u$ is given by a Barenblatt solution.

Consider the
rescaled function
\begin{equation}\label{eqn-tilu}
\tilde{u}(x,\tau) =
(T-t)^{-\beta} u(x\, (T-t)^{\gamma},t), \qquad \tau = - \log (T-t)
\end{equation}
with $\beta$ and $\gamma$ given by \eqref{eqn-parameters}. It follows
by direct computation, that $u$
satisfies  the equation
\begin{equation}
\label{eqn-rescaled}
\tilde{u}_\tau = \Delta \tilde{u}^m + |\gamma| \, \dv (x\cdot \tilde{u})
\end{equation}
and due to condition (\ref{eqn-trapped}), the inequality
\begin{equation}
\label{eqn-rescaled-trap}
\left (\frac{C^*}{|x|^2 + k_1} \right )^{\frac{1}{1-m}} \le \tilde{u}(x,\tau) \le
\left (\frac{C^*}{|x|^2 + k_2} \right )^{\frac{1}{1-m}}
\end{equation}
holds, for $(x,\tau) \in \R^N\times [-\log T,\infty)$.
We denote by
\begin{equation}
\label{eqn-rbaren}
\tilde{B}_{k}(x) = \left ( \frac{C^*}{k+ |x|^2} \right )^{\frac{1}{1-m}}
\end{equation}
the rescaled Barenblatt solution.

Our convergence results are described in the  following two
Theorems. The first result is concerned with the range of
exponents $\frac{N-4}{N-2}< m < \frac{N-2}{N}$,  for which the
difference $u-B_k$, with $B_k$ a Barenblatt solution,   is
integrable, amely we have  $\int_{\R^N} (u-B_k) (x,t) \, dx <
\infty$. Notice that this range of exponents includes the  Yamabe
flow, $m=(N-2)/(N+2)$, when $N < 6$.

\begin{thm}
\label{thm-integrable}
Let $u$ solve the equation  \eqref{eqn-u}  for $
\frac{N-4}{N-2} < m < \frac{N-2}{N}$, with initial value $u_0$ satisfying
(\ref{eqn-trapped}), for some constants $k_1, k_2$. Then,
the rescaled function $\tu$ given by \eqref{eqn-tilu}
converges, as $\tau\to\infty$,  uniformly on $\R^N$, and also in $L^1(\mathrm{R}^N)$,
to the rescaled Barenblatt solution $\tilde B_{k_0}$ given by
\eqref{eqn-rbaren}, for some $k_0 >0$.
The constant $k_0$ is uniquely determined by the equality
$$\int_{\R^N} u_0 \, dx = \int_{\R^N}  B_{k_0}  \, dx.$$
\end{thm}

 The second result deals with the
 range of exponents $0< m
\leq \frac{N-4}{N-2}$, for which the difference $u-B_k$, with $B_k$ a
Barenblatt solution, is non-integrable, namely  $\int_{\R^N} (u-B_k)
(x,t) \, dx = \infty$.

\begin{thm}
\label{thm-nonint}
Let $u$ solve the equation \eqref{eqn-u}, for $0 < m \le
\frac{N-4}{N-2}$, $N >4$, with initial value $u_0$ satisfying
(\ref{eqn-trapped}), for some constants $k_1, k_2$. Assume, in addition that
\begin{equation}
\label{eqn-int-diff}
u_0 = B_{k_0} + f
\end{equation}
for some $k_0 >0$, where $B_{k_0}$ is a Barenblatt solution and
$f$ is in $L^1({\R}^N)$. Then, the rescaled function $\tu$ given
by \eqref{eqn-tilu} converges, as $\tau\to\infty$, uniformly on
$\mathrm{R}^N$, to the Barenblatt solution $\tilde B_{k_0}$.
\end{thm}

\begin{rem}
\label{remark-integrable}
The condition $0 < m \le \frac{N-4}{N-2}$
in Theorem \ref{thm-nonint} implies that for any two Barenblatt solutions
$B_k$ and $B_{k'}$, $B_k - B_{k'} \in L^1(\R^N)$ if and only if $k = k'$.
\end{rem}

One may ask whether condition \eqref{eqn-igc1} is necessary for
Theorems \ref{thm-integrable} and \ref{thm-nonint} to hold true. We
will show in section \ref{sec-example} that this is indeed the case.
We will present an example of a class of initial conditions $u_0$
which satisfy the growth condition
\begin{equation}\label{eqn-igc2}
u_0(x) = \left ( \frac{C^* \, T}{|x|^2}  \right )^{\frac 1{1-m}} (1 + o(1)), \qquad \mbox{as} \,\, |x| \to \infty
\end{equation}
with  $C^*$  given by \eqref{eqn-C*},  for which the solution $u$
of  \eqref{eqn-u}  with initial data $u_0$ satisfies the following Theorem.

\begin{thm}
\label{thm-example}
There exists a class of  solutions  $u$  of  the Cauchy problem  \eqref{eqn-u}
with initial data $u_0$ satisfying \eqref{eqn-igc2} and with the following
properties:
\begin{enumerate}
\item[(i)]
The vanishing time $T^*$ of $u$ satisfies  $T^* > T$.
\item[(ii)]
The solution $u$  satisfies as $|x| \to \infty$,  the growth conditions
\begin{equation}\label{eqn-igc3}
u(x,t) \geq  \left (  \frac{C^* \, (T-t) }{1+|x|^2} \right )^{\frac 2{1-m}}, \qquad   \mbox{on} \,\,\,  0 < t < T \end{equation}
and
\begin{equation}\label{eqn-igc4}
u(x,t) \leq \frac{C(t)}{|x|^{\frac m{N-2}}}, \qquad \mbox{on} \,\,
T < t <T^*. \end{equation}
In particular, $u$ becomes integrable on  $t >T$.
\item[(iii)] The vanishing behavior of $u$ is given by one of the self-similar solutions $\Theta(x,t)$ (see
section \ref{sec-example} for the explanation of $\theta(x,t)$).
\end{enumerate}
\end{thm}

The vanishing behavior of the solution $u$ in this case is described in the results of Galaktionov and Peletier \cite{GPe},
and del Pino and S\'aez \cite{DS} (in the case $m =  \frac{N-2}{N+2}$), also formally shown by King \cite{King}.

The case of a special interest is $m = \frac{N-2}{N+2}$, when the
equation for $u$ is eqivalent to the Yamabe flow of a corresponding
conformally flat metric. The previous Theorem gives the following
corrolary in the Yamabe case.

\begin{cor}
\label{cor-yamabe}
If $u$ is a solution of the Cauchy problem \eqref{eqn-u} with initial
data satisfying \eqref{eqn-igc2}. Then the vanishing time of $u$ is $T^* >
T$, (ii) in Theorem \ref{thm-example} holds and
$$(T^* - t)^{-\frac{1}{1-m}}u(x,t) \to (\frac{C_1}{C_2 + |x-\bar{x}|^2})^{\frac{N+2}{2}},$$
as $t\to T^*$, where $\bar{x}\in \R^N$ and $C_1, C_2 > 0$.
\end{cor}

Geometrically speaking, Yamabe flow starting at $u_0$ (described
above) developes a singularity at $T < T^*$ at which the Barenblatt
solution (cylinder) pinches off and immediately at $t > T$ the
solution becomes integrable. Due to the results of Del Pino and Saez
it smoothens up at $t > T$ and exhibits the behaviour of a compact
sphere as $t\to T^*$.

\smallskip

The next section will be devoted on preliminary estimates for
solutions  $u$ of \eqref{eqn-u} with initial data in
$L^1_{\mbox{loc}}$. The proof  of Theorem \ref{thm-integrable}
will be given in section \ref{sec-integrable}. It will follow from
the  strong  $L^1$-contraction principle, Lemma
\ref{lem-contraction}, which holds for  the difference of the
rescaled solutions $\tilde u - \tilde B_k$, for any   Barenblatt
solution $B_k$.  This method, based upon  the ideas of Osher and
Ralston \cite{OR}, was previously  used  by S-Y Hsu in \cite{Hs1}.
Since the difference $\tilde u - \tilde B_k \notin L^1(\R^N)$  in
the range of exponents $0 < m < \frac {N-4}{N-2}$, for the proof
of Theorem  \ref{thm-nonint} we will need to weight the
$L^1$-norms with an appropriate power $\tilde B_{k_2}^\alpha$. The
proof in this other case is more involved and will be given in
section \ref{sec-noint}. The last section will be devoted to the
construction of the examples described in Theorem \ref{thm-example}.

{\bf Acknowledgements:} The authors would like to thank  S. Brendle,
M. Gr\"uneberg, J.R. King and J.L. Vazquez  for helpful discussions. The second author would
like to thank Max Planck in Golm (Potsdam) for the hospitality during
the period in which part of this work has been done.

\end{section}

\begin{section}{Preliminary Estimates}\label{sec-prelim}

Our goal in this section is to establish  the  $L^1-$contaction
for solutions of \eqref{eqn-u} and \eqref{eqn-rescaled} and some
other preliminary results which can also be of independent
interest. We begin by showing the following integrability lemma
for the difference of any two solutions of \eqref{eqn-u}.

\begin{lem}
\label{lem-integrability}
Assume that $u, v$ are two solutions of \eqref{eqn-u}
on $\R^N \times (0,T)$.
If $f=u_0- v_0 \in L^1(\R^N)$ and $f$ is compactly supported, then $u(\cdot,t) - v(\cdot,t) \in L^1(\R^N)$
for all $t\in [0,T)$.
\end{lem}

\begin{proof} Our proof is based on the well known technique of Herrero and Pierre
\cite{HP}.  We introduce  the potential function
 \begin{equation}\label{eqn-potential}
 w(x,t) = \int_0^t | (u^m  - v^m)(x,s)  | \, ds
 \end{equation}
which  satisfies the inequality
\begin{equation}
\label{equation-subsolution}
\Delta w \geq - |f|, \qquad \mbox{on}
\,\,\, \R^N.
\end{equation}
Indeed,  by Kato's inequality
\cite{K}, we have
$$\Delta |u^m - v^m | \ge \sign \, (u-v)\, \Delta(u^m - v^m)$$
so from equation   (\ref{eqn-u}), we obtain
\begin{equation}
\label{eqn-abs-value}
\frac{\partial}{\partial t}\,  |u - v| \le \Delta |u^m - v^m|.
\end{equation}
Integrating  the previous inequality in time, and using that  $|f|
= |u_0 -v_0|$, we obtain \eqref{equation-subsolution}.

Let
$$Z(x) = \frac 1{N(N-2)\omega_N} \int_{\R^N} \frac 1{|x-y|^{N-2}} \, |f(y)| \, dy$$
denote  the Newtonian potential of $|f|$, so that from  \eqref{equation-subsolution} we have
$$\Delta (w - Z) \ge 0.$$
Also, since $|f|$ is integrable and compactly supported, there exists
a constant $C < \infty$ for which
\begin{equation}\label{eqn-estZ}
Z(x) \leq \frac{C}{|x|^{N-2}}, \qquad \forall x >> 1,
\end{equation}
where the range of $x$ for which the estimate holds depends on the
support of $f$. The mean value property implies that
\begin{equation}
\begin{split}
\label{eqn-MVP}
w(x) &\le Z(x) + \frac{1}{\omega_N \rho^N}\int_{B_\rho(x)} (w(y) - Z(y))\, dy \nonumber \\&\le  Z(x) +  \frac{1}{\omega_n \rho^N}\int_{B_\rho(x)}w(y)\,dy
\end{split}
\end{equation}
for all $x \in \R^N$, $\rho >0$.

We next claim that
\begin{equation}
\lim_{\rho\to\infty}\frac{1}{\rho^N}\int_{B_{\rho(x)}}w(y)\,dy = 0.
\end{equation}
Indeed, by Lemma $3.1$ in \cite{HP} we have
$$\int_{B_{\rho}(x)} |u - v|(y,t)\, dy \le
C\, \left ( \| \, u_0 - v_0 \, \|_{L^1(\R^N)} + \rho^N \, \left (\frac t{\rho^2} \right )^{\frac{1}{1-m}} \right)$$
which yields
\begin{eqnarray}
\label{eqn-better-1}
\int_{B_{\rho}(x)}w(y,t) \, dx &=&  \int_0^t\int_{B_{\rho}(x)} |u- v|^m(y,s)\,  \, dy \, ds \nonumber \\
&\le& C\int_0^t \rho^{N(1-m)} \left ( \int_{B_{\rho}(x)} |u - v| \, dy \right )^m \, ds \nonumber \\
&\le& C\rho^{N(1-m)}\int_0^t \left ( \| f \|_{L^1(\R^N)} + \rho ^N \, \left (\frac  s{\rho^2} \right )^{\frac{1}{1-m}} \right)^m \, ds   \nonumber \\
&\le& C(T, \|f\|) \, \rho^{N-2m/(1-m)}.
\end{eqnarray}
Combining \eqref {eqn-estZ}, (\ref{eqn-MVP}), (\ref{eqn-better-1}) and letting $\rho\to\infty$, we conclude the estimate
\begin{equation}
\label{eqn-decay}
w(x) \le Z(x) \le \frac{C}{|x|^{N-2}}.
\end{equation}

We will use \eqref{eqn-decay} along the lines of the proof of Theorem 2.3  in \cite{HP} to bound $\| (u-v)(\cdot,t)\|_{L^1(\R^N)}.$
For $0 \le \eta_R \le 1$, let $\eta_R \in C_0^{\infty}(\R^N)$ be a test function such that $\eta_R = 1$
for $|x| \le R$  and $\eta_R = 0$ for $|x| \ge 2R$. Then $|\Delta\eta_R| \le C/R^2$ and
$|\nabla\eta_R| \le C/R$. Using  equation \eqref{eqn-u},  estimate  (\ref{eqn-decay}) and the integrability of $f$, we get
\begin{eqnarray*}
\int |u-v|(\cdot ,t)\, \eta_R \,dx &\le& \int_{\R^N} |f|\, dx +
\int_0^t\int_{B_{2R} \setminus  B_R} |u^m - v^m| \, \Delta \eta_R \,  dx\, ds \\
&\le& \|f\|_{L^1(\R^N)} + \frac{C}{R^2}\int_{B_{2R}\backslash B_R} w(x,t)\, dx \\
&\le& C(\|f\|_{L^1(\R^N)}) + \bar{C}
\end{eqnarray*}
where $\bar{C}$ can be taken to be independent of $R$ because of \eqref{eqn-decay}.
Letting $R\to \infty$ in the previous estimate gives that
$$\sup_{t\in [0,T)}\int_{\R^N} |u-v|(x,t)\, dx \le C(\|f\|_{L^1(\R^N)}) + \bar{C} < \infty$$
finishing the proof of the lemma.
\end{proof}

As a   consequence of the previous Lemma, we will now establish
the following $L^1$ contraction principle for the solutions to
(\ref{eqn-u}) that are bounded from below by a Barenblatt solution
$B$.

\begin{cor}
\label{cor-contraction}
Let  $u$, $v$ be  two solutions of \eqref{eqn-u} with initial values $u_0$,
$v_0$ respectively and that $f=u_0 - v_0 \in L^1(\R^N)$. Assume in addition that  $u, v \ge B$, for some Barenblatt solution $B$ given by \eqref{eqn-baren}. Then,
\begin{equation}\label{eqn-contract}
\int_{\R^N}|u(\cdot,t) - v(\cdot,t)|\, dx \le \int_{\R^N}|u_0 - v_0|\, dx,
\qquad \forall  t\in [0,T).
\end{equation}
\end{cor}

\begin{proof}
Let $\eta_R$ be a cut off function as in the proof of Lemma
\ref{lem-integrability}, with the support contained in $B_{2R}$. Then,
as before,  we have
$$\frac{\partial}{\partial t}|u(x,t) - v(x,t)| \le \Delta|u^m - v^m|.$$
If we multiply the above  inequality by  $\eta_R$ and integrate over $\R^N$, since
$|\Delta \eta_R| \le \frac{C}{R^2}$, we get
$$\frac{d}{dt}\int_{\R^N}|(u - v)(x,t)| \, \eta_R\, dx \le \frac{C}{R^2}
\int_{B_{2R}\backslash B_R}a(x,t)\, |(u - v)(x,t)|$$ with
$$a(x,t) = \int_0^1\frac{d\theta}{(\theta u + (1-\theta)v)^{1-m}} \le C R^2, \qquad \mbox{in} \,\, B_{2R}$$
 (since $(\theta u + (1-\theta)v)^{1-m} \ge B^{1-m}
= C_1(|x|^2 + C_2)$).
Hence, fixing   $t\in [0, T)$, we obtain the estimate
$$\frac{d}{dt}\int_{\R^N}|u - v|(x,t)\eta_R\, dx  \le
C\int_{B_{2R}\backslash B_R}|u - v|(x,t)\, dx.$$
Assume that $u_0, v_0$ are compactly supported, so that $f=u_0-v_0$ is compactly supported as well. Then,  the right hand side of the above inequality converges to zero as $R\to\infty$, due to
Lemma \ref{lem-integrability}. This gives
$$\frac{d}{dt}\int_{\R^N}|u - v|(x,t)\, dx \le 0$$
which implies \eqref{eqn-contract}.

To remove the assumption that $u_0$, $v_0$ are compactly
supported, we use a standard approximation argument. For any $k
>1$, we set  $u_0^k=u_0 \, \chi_{B_k(0)}$, $v_0^k=v_0 \,
\chi_{B_k(0)}$. Let  $u^k$, $v^k$ be the solutions  of
\eqref{eqn-u} on $\R^N \times (0, \infty)$ with initial values
$u_0^k$, $v_0^k$ respectively. By standard arguments $u^k \to u$
and $v^k \to v$ uniformly on compact subsets of $\R^n \times
(0,\infty)$. Also,
 by the previous argument
$$\int_{\R^N}|u^k(\cdot,t) - v^k(\cdot,t)|\, dx \le \int_{\R^N}|u_0^k - v_0^k|\, dx \leq \int_{\R^N}|u_0 - v_0|\, dx$$
for all $ t >0$ and for all $k$. Letting $k \to \infty$ we readily obtain \eqref{eqn-contract}.
\end{proof}

As an immediate corollary of Lemma \ref{lem-integrability} we have the following result concerning the rescaled solutions $\tilde u$, $\tilde v$.

\begin{cor}
\label{cor-integrability-rescaled}
Let $u$ and $v$ be as above. If $u_0 - v_0\in L^1(\R^N)$,  then
for every $\tau > -\log T$ there is $C(\tau)$ such that
$$\int_{\R^N}|\tilde{u}(x,\tau) - \tilde{v}(x,\tau)|\,dx \le C(\tau).$$
\end{cor}

Let $u, v$ be two  solutions of \eqref{eqn-u} satisfying \eqref{eqn-trapped}.  We
observe that the difference  $q = \tilde{u} - \tilde{v}$ satisfies the equation
\begin{equation}
\label{eqn-rescale-v}
q_{\tau} = \Delta(a(x,\tau)\, q) + |\gamma| \, \dv(x \cdot q(x,\tau))
\end{equation}
on $ \R^N\times [-\log T,\infty)$, with
\begin{equation}\label{eqn-a}
a(x,\tau) = \int_0^1\frac{m}{(\theta \, \tilde{u} + (1-\theta)\, \tilde{v})^{1-m}}\, d\theta.
\end{equation}
Since both  $\tilde{u} $ and $ \tilde{v}$ satisfy (\ref{eqn-rescaled-trap}),  for some
constants $k_1, k_2$, it is clear that $a(x,\tau)$ is smooth  and satisfies the growth estimate
\begin{equation}\label{eqn-a2}
\frac{m\, (k_2+|x|^2)}{C^*} \le a(x,\tau) \le \frac{m\, (k_1+|x|^2)}{C^*}.
\end{equation}
Hence, (\ref{eqn-rescale-v}) is uniformly parabolic on any compact subset of $\R^N \times (0,\infty)$.

Let   $F(x,t)$  be a solution of
\begin{equation}
\label{eqn-heat}
F_t = \Delta (a_1(x,t)\, F)
\end{equation}
with
$$a_1(x,t) = \int_0^1\frac{m}{(\theta \, u(x,t) + (1-\theta)\, v(x,t))^{1-m}}\, d\theta.$$
A direct computation shows that   $\tilde{F}(x,\tau)= F(x,(T-t)^{\gamma},t)$,   $\tau = -\log(T-t)$,  is a solution of equation
$$\tilde{F}_{\tau} = \Delta(a(x,\tau)\, \tilde{F}) + |\gamma| \, \dv (\tilde{F}\cdot x)$$
where $a(x,\tau)$ is given by \eqref{eqn-a}. Similarly as before we have the following result.

\begin{cor}
\label{cor-heat-int}
Let $F(x,t)$ be  a solution of (\ref{eqn-heat}). If $F(x,0)\in L^1(\R^N)$,  then
$F(x,t) \in L^1(\R^N)$. Moreover, $\tilde{F}(x,\tau)\in L^1(\R^N)$ and
for every $\tau > - \log T$ there is a $C(\tau)$ such that
$$\|F(\cdot,\tau)\|_{L^1(\R^N)} \le C(\tau).$$
\end{cor}
\end{section}

%We finish this section with the following result which provides a lower bound on the vanishing time
%of solutions of \eqref{eqn-u}  with initial data satisfying the growth condition \eqref{eqn-igc1}.

%\begin{prop} Let $u$ be a solution of \eqref{eqn-u} with initial data $ u_0  > 0$ satisfying the growth
%condition
%\begin{equation}\label{eqn-igc100}
%u_0(x) = \left ( \frac{C^*\, T}{|x|^2}  \right )^{\frac 1{1-m}} (1 + o(1)), \qquad \mbox{as} \,\, |x| \to \infty.
%\end{equation}
%with $C^*$ given by \eqref{eqn-C*}. Then, the vanishing time $T^*$ of $u$ satisfies $T^* \geq T$.
%\end{prop}

%
%\begin{proof}
%Assume first that $u_0 >0$ on $\R^N$Let $\delta >0$ be a small number. Since $u_0$ satisfies \eqref{eqn-igc100}, there exists $\rho_{\delta}  >0$ for which $u_0(x) \geq ({C^*\, (T-\delta)}/{|x|^2}  )^{\frac 1{1-m}}

%\end{proof}
%
%
\begin{section}{The integrable case}\label{sec-integrable}

This section is devoted to the proof of Theorem
\ref{thm-integrable} which deals with solutions of equation
\eqref{eqn-u} in the range of exponents $\frac{N-4}{N-2} < m <
\frac{N-2}{N}$. In this case the difference of two solutions $u,v$
satisfying \eqref{eqn-trapped2} is integrable. We begin this
section with the following strong contraction principle, which
constitutes the main step in the proof of Theorem
\ref{thm-integrable}. Its proof as well as the rest of the
argument is very similar to the proof of Theorem $2.3$ in
\cite{Hs1}. To facilitate future references  we will sketch the proof of
the strong contraction principle.

\begin{lem}
\label{lem-contraction2}
Let $u,v$ be two solutions of \eqref{eqn-u}, for $m \in
(0,\frac{N-2}{N}$, with initial values $u_0,v_0$, satisfying
(\ref{eqn-trapped}). If
$$\min\{\|(\tilde{u}_0 - \tilde{v}_0)_+\|_{L^{\infty}}, \|(\tilde{v}_0 - \tilde{u}_0)_-\|_{L^{\infty}}\} > 0,$$
then
$$\|(\tilde{u} - \tilde{v})(\cdot,\tau)\|_{L^1(\R^N)} < \|\tilde{u}_0
- \tilde{v}_0\|_{L^1(\R^N)}, \qquad \tau \ge -\log T.$$
\end{lem}

\begin{proof}
Notice that by the comparison principle, $\tilde{u}(x,\tau)$,
$\tilde v(x,\tau)$ satisfy (\ref{eqn-rescaled-trap}). The proof is
almost the same as that of Lemma $2.1$ by S-Y Hsu in \cite{Hs1},
using the results in  \cite{OR} and \cite{Zh}.

Set $q=\tilde u - \tilde v$ and observe, as above, that $q$
satisfies equation \eqref{eqn-rescale-v}.  Fix $R >0$. By the
standard parabolic theory, there exist  solutions $q_+^R, q_-^R$
of \eqref{eqn-rescale-v} in $Q_R = B_{R}\times (-\log T,\infty)$,
with initial values $q(\cdot,-\log T)_+, \, q(\cdot,-\log T)_-$
and boundary values $q_+$, $q_-$ on $\partial B_R \times
(0,\infty)$, respectively. Notice that $q_+^R - q_-^R$ is a
solution of (\ref{eqn-rescale-v}) in $Q_R$, with initial value
$q(\cdot,-\log T)$ and boundary values $q_+ - q_-$. By the maximum
principle, $q = q_+^R - q_-^R$ on $Q_R$. Similarly there are
solutions $\bar{q}_+^R$ and $\bar{q}_-^R$ of (\ref{eqn-rescale-v})
in $Q_R$ with initial values $q_+$, $q_-$, and zero lateral
boundary value. By the maximum principle, $0 \le \bar{q}_+^R \le
q_+^R$ and $0 \le \bar{q}_-^R \le q_-^R$. Furthermore let
$\tilde{q}_R$ be the solution of (\ref{eqn-rescale-v}) with
initial value and lateral boundary value $\tilde{B}_{k_2} -
\tilde{B}_{k_1}$. By the maximum principle, we have
$$0 \le q_+^R, q_-^R \le \tilde{q}_R.$$
Let $\eta \in C_0^{\infty}(\R^2)$ be a cut-off function
such that $\eta(x) = 1$ on  $|x| \le 1/2$, $\eta(x) = 0$ for all $|x|
\ge 1$ and $0 \le \eta \le 1$.  Denote by $\eta_R = \eta(x/R)$. The
same computation as in the proof of Lemma $2.1$ in \cite{Hs1} gives
\begin{eqnarray*}
& &\int_{\R^N}|q(x,\tau)|\, \eta_R\,dx - \int_{\R^N} |\tilde{u}_0 - \tilde{v}|\, \eta_R \,dx\nonumber \\
&=& \int_{-\log T}^{\tau}\left (\int_{R^N} a(x,\tau')\, q_+^R(x,\tau')\, \Delta \eta_R  -
|\gamma| \, q_+^R(x,\tau')\, x\cdot \nabla\eta_R\,dx \right ) d\tau'\nonumber \\
&+& \int_{-\log T}^{\tau} \left (\int_{R^N}a(x,\tau')\, q_-^R(x,\tau')\, \Delta \eta_R -
|\gamma| \, q_-^R(x,\tau') \, x\cdot\nabla\eta_R\,dx \right ) d\tau' \nonumber\\
&& \,\, -  2\int_{\R^N}\min\{q_+^R(x,\tau'), q_-^R(x,\tau')\}\, \eta_R\, dx.
\end{eqnarray*}
The families of solutions   $\bar{q}_+^R(x,\tau)$ and $\bar{q}_-^R(x,\tau)$ are
monotone increasing in $R$ and  uniformly bounded above, which
implies that
$$\bar{q}_1 = \lim_{R\to\infty}\bar{q}_+^R$$
and
$$\bar{q}_2 = \lim_{R\to\infty}\bar{q}_-^R$$
exist and are both solutions of (\ref{eqn-rescale-v}) on $\R^N\times (-\log T,\infty)$.
This implies
\begin{eqnarray}
\label{eqn-limit}
& &\int_{\R^N}|q(x,\tau)|\eta_R\,dx - \int_{\R^N}|\tilde{u}_0 - \tilde{v}_0|\eta_R \,dx \nonumber \\
&\le& \frac{C}{R^2}\int_{-\log T}^{\tau}\int_{R/2\le|x|\le R}\, a(x,\tau')\tilde{q}^R(x,\tau')\,dx d\tau' \nonumber \\
&+&  C\int_{-\log T}^{\tau}\int_{R/2\le|x|\le R}\tilde{q}^R(x,\tau')\,dx d\tau' - 2\int_{B_{R_0}}\min\{q_+^R(x,\tau), q_-^R(x,\tau)\}\eta_R\, dx. \nonumber \\
&\le& C\int_{-\log T}^{\tau}\int_{R/2\le|x|\le
R}\tilde{q}^R(x,\tau')\,dx d\tau' -
2\int_{B_{R_0}}\min\{\bar{q}_+^R(x,\tau),
\bar{q}_-^R(x,\tau)\}\eta_R\, dx.\nonumber
\end{eqnarray}
By the same computation as in \cite{Hs1}, after letting
$R\to\infty$,  we get
\begin{equation*}
\begin{split}
\int_{\R^N}|q(x,\tau)|\,dx &- \int_{\R^N}|\tilde{u}_0 - \tilde{v}_0| \,dx \\
&\le - 2\int_{B_{R_0}}\min\{\bar{q}_1(x,\tau),
\bar{q}_2(x,\tau)\}\, dx
\end{split}
\end{equation*}
for all $ R_0 >0$ and $ \tau> -\log T$.
Since $\bar{q}_1 \ge q_+^{2R_0}$ and $\bar{q}_2 \ge q_-^{2R_0}$, we
obtain
\begin{equation*}
\begin{split}
\int_{\R^N}|q(x,\tau)|\,dx &- \int_{\R^N}|\tilde{u}_0 -\tilde{v}| \,dx\\
&\le - 2\int_{B_{R_0}}\min\{q_+^{2R_0}(x,\tau),
q_-^{2R_0}(x,\tau) \}\, \, dx.
\end{split}
\end{equation*}
Since  $q_+^{2R_0}$ and $q_-^{2R_0}$ are the solutions of (\ref{eqn-rescale-v}) in
$Q_{2R_0}$ with zero boundary value and initial values $q_+(\cdot,-\log T)$,
$q_-(\cdot,-\log T)$, respectively, by the Green's  function representation for solutions, for any $\tau > -\log T$,  there exists a constant $c(\tau)$ such that
$$\min_{x\in B_{2R_0}}q_+^{2R_0} \ge c(\tau) > 0 \quad \mbox{and} \quad \min_{x\in B_{2R_0}}q_-^{2R_0} \ge c(\tau) > 0. $$
That finishes the proof of Lemma \ref{lem-contraction}.
\end{proof}

The rest of the proof relies on the following result of Osher and
Ralston (\cite{OR}).

\begin{lem}[Lemma 1 in \cite{OR}]
\label{lem-or1}
Suppose that $\tilde{u}(\cdot,\tau_i)  \to \bar u$ in $L^1(\R^N)$,
as $i \to \infty$, for some sequence $\tau_i \to \infty$. Let   $\tilde{B}_k$
be any stationary solution of (\ref{eqn-rescaled}). If $\tilde{v}$
is the  solution of (\ref{eqn-rescaled}) in $\R^N\times [0,\infty)$
with initial value $\tilde{v}(x,0) = \bar u(x)$, then
$$\|\tilde{v}(\cdot,\tau) - \tilde{B}_k\|_{L^1(\R^N)} =
\|\bar u  - \tilde{B}_k\|_{L^1(\R^N)}, \quad \forall \tau >0, \,\, k >0.$$
\end{lem}

\noindent{\em Proof of Theorem \ref{thm-integrable}}.
We claim that in this integrable
case there is a unique $k_0$ so that
$$\int_{\R^N} (u_0 - B_{k_0})\, dx = 0.$$
To prove that, let
$$f(k) = \int_{\R^N}(u_0 - B_k)\, dx$$
and observe that $f(k)$ is a continuous, monotone increasing function
with $f(k_1) \ge 0$ and $f(k_2) \le 0$ due to (\ref{eqn-trapped}).
%Fix $k'$ and write
%$$f(k) = f(k') + \int_{\R^N}(B_{k'} - B_k)\, dx.$$
%As $k\to\infty$, $B_{k'} - B_k > 0$ and it is monotone
%increasing in $k$. Therefore, by the Lebesgue
%monotone convergence theorem,
%\begin{eqnarray*}
%\lim_{k\to\infty}f(k) &=& f(k') + \int_{\R^N}(B_{k'} - \lim_{k\to\infty}B_k)\,dx \\
%&=& f(k') + \int_{R^N}B_{k'}\,dx = \infty.
%\end{eqnarray*}
%Similary, as $k\to 0$, $B_k - B_{k'} > 0$ and it is monotone increasing, so by the Lebesgue monotone convergence theorem,
%\begin{eqnarray*}
%\lim_{k\to 0}f(k) &=& f(k') - \int_{R^N}(\lim_{k\to 0}B_k - B_{k'})\,dx  \\
%=& \int_{R^N} u_0 -  \left ( \frac{C^*}{|x|^2} \right )^{\frac 1{1-m}} \,dx \\
%&\le& \int_{R^N} \left (\frac{C^*}{|x|^2 + k_2} \right )^{\frac 1{1-m}} -
%&\left ( \frac{C^*}{|x|^2} \right )^{\frac 1{1-m}} \,dx < 0.
%\end{eqnarray*}
Therefore, by the intermediate value theorem there exists a unique
$k_0$ such that $ f(k_0) = 0.$
The rest of the
proof is the same as in \cite{Hs1}, based on the strong contraction  Lemma \ref{lem-contraction2}.
\end{section}

\begin{section}{The non-integrable case}
\label{sec-noint}

This section will be devoted to the proof of Theorem
\ref{thm-nonint} which is concerned with the range of exponents $0
< m \leq \frac{N-4}{N-2}$, $N >4$. Assume that $u$ is a solution
of \eqref{eqn-u} which satisfies the bound \eqref{eqn-trapped2}.
Throughout this section $\tilde u$ will denote the rescaled
solution defined by \eqref{eqn-tilu},  and $\tilde B_k$ the
rescaled Barenblatt solution given by \eqref{eqn-rbaren}.

Since a difference of any two solutions $ u$, $v$ of equations
\eqref{eqn-u} is not always integrable in the range of exponents
$0 < m \leq \frac{N-4}{N-2}$, $N >4$, we need to depart in this
section from the techniques used in \cite{Hs1} in which we heavily
used the integrability of a difference of any two Barenblatt
solutions.

We
define the  weighted $L^1$-space with  weight $\tilde{B}^{\alpha}:=
(C^*/(k_2+|x|^2))^{\alpha/(1-m)}$, as
$$L^1(\tilde{B}^{\alpha}, \R^N)
:= \{f | \int_{\R^N} f (x) \, \tilde{B}^{\alpha}(x)\,dx < \infty\}.$$

%Similarly to the definition above,  we define the space $$L^1(\tilde{B}^{\alpha} \, \eta_R, B_R):= \{f | \int_{\R^N} f \, \tilde{B}^{\alpha}(x)\,dx < \infty\}.$$

\begin{lem}
\label{lem-contraction}
Let $u,v$ be two solutions of (\ref{eqn-rescaled}),
 for $0 < m <
\frac{N-4}{N-2}$, with initial values $u_0, v_0$ satisfying
(\ref{eqn-trapped}) and $u_0 - v_0 \in L^1(\tilde{B}^{\alpha},
\R^N)$, with  $\alpha = (N-2)(1-m)/2-1$ (we will  explain later
such a choice of $\alpha$). Let $\tilde{B} := \tilde{B}_{k_2}$.
If
\begin{equation}
\label{equation-min-max}
\max_{\R^N}|\tilde{u}_0 - \tilde{v}_0| \neq 0
\end{equation}
then,
$$\|(\tilde{u} - \tilde{v})(\cdot,\tau) \, \tilde{B^{\alpha}}\, \eta_R\|_{L^1} < 
\|(\tilde{u}_0- \tilde{v}_0)\, \tilde{B}^{\alpha}\eta_R\|_{L^1}, \quad \tau \ge -\log T.$$
\end{lem}

\begin{proof}
Set $q=\tilde u - \tilde v$. After rescaling,
(\ref{eqn-abs-value}) becomes
$$|q|_{\tau} \le \Delta(a|q|) + |\gamma|\nabla(x\cdot|q|),$$
where $a(x,\tau)$ is as in (\ref{eqn-a}). Let $\eta \in
C_0^{\infty}(\R^2)$ be a cut off function as before. Denote by
$\eta_R = \eta(x/R)$ so that $|\nabla \eta_R| \leq C/R$, $|\Delta
\eta_R| \leq C/R^2$.
Then, the above equation and integration by parts    yield to
\begin{equation*}
\begin{split}
 \int_{\R^N} &|q(x,\tau) |\eta_R\,  \tilde{B}^{\alpha}(x)\,dx -
\int_{\R^N}|\tilde{u}_0 - \tilde{v}_0 | \, \eta_R \, \tilde{B}^{\alpha}(x) \,dx =\\
&= \int_{-\log T}^{\tau} \int_{R^N}  \left ( a(x,\tau')\,
|q|(x,\tau') \, \{\tilde{B}^{\alpha}(x) \, \Delta \eta_R
+ \Delta\tilde{B}^{\alpha}(x)\, \eta_R + 2 \nabla\eta_R \, \nabla \tilde{B}^{\alpha}(x) \} \right .  \nonumber \\
&- \left .  |\gamma| \, |q|(x,\tau')\, x\cdot \, \{
\tilde{B}^{\alpha}(x)\, \nabla\eta_R+ \nabla\tilde{B}^{\alpha}(x)\,
\eta_R \}\,dx \right )  d\tau'. \nonumber \\
\end{split}
\end{equation*}
Moreover,
\begin{equation}
\begin{split}
\label{equation-bounded-q}
\quad \int_{\R^N}|q(x,\tau)| &\, \eta_R
\tilde{B}^{\alpha}(x)\,dx
- \int_{\R^N}|\tilde{u}_0 - \tilde{v}_0| \, \eta_R\, \tilde{B}^{\alpha}(x) \,dx \\
&\le \frac{C}{R^2}\int_{-\log T}^{\tau}\int_{R/2 \le |x| \le
R}a(x,\tau')|q|(x,\tau')\tilde{B}^{\alpha}(x)\,dx d\tau' \\
&+ \frac{C}{R }\int_{-\log T}^{\tau}\int_{R/2 \le |x| \le R}
a(x,\tau') |q|(x,\tau')\, |\nabla\tilde{B}^{\alpha}(x)| \,dx d\tau'  \\
&+ C\int_{-\log T}^{\tau}\int_{R/2 \le |x| \le
R}|q|(x,\tau')\tilde{B}^{\alpha}(x)\,dx d\tau'  \\
&+ \int_{-\log T}^{\tau}\int_{\R^N} \{
a(x,\tau')\Delta\tilde{B}^{\alpha} -
|\gamma| \,  x \cdot \nabla\tilde{B}^{\alpha} \} \,|q|(x,\tau')\,\eta_R \,  dx d\tau' \\
&= I_1 (R) + I_2 (R) + I_3 (R) + I_4 (R)\qquad \forall \,\, R\ge R_0 > 0. \\
\end{split}
\end{equation}
We fix
\begin{equation}
\label{equation-value} \alpha = \frac{(N -2)\,(1- m)}{2}-1.
\end{equation}

\noindent{\em Claim i.}  If $\max_{\R^N}|\tilde{u}(x,\tau) - \tilde{v}(x,\tau)| \neq 0$,
there exists a constant $C(\tau) >0$, such that $I_4(R)  < -C(\tau) $.

\noindent{\em Proof of Claim i.} We recall that
$$\frac{m\, (|x|^2 + k_2)}{C^*} \le a(x,\tau) \le \frac{m\, (|x|^2 + k_1)}{C^*}$$
with $C^*=\frac{2m\, (N - 2 -mN)}{1-m}.$
A direct computation shows that
$$\Delta B^\alpha = \frac{2\, \alpha \, (\, [ 2\alpha  - (1-m)\, (N-2)]  \, |x|^2 - k_2\, (1-m)\, N) }{(1-m)^2 \, (1+ |x|^2)^2} \, B^\alpha \leq 0$$
provided that $ \alpha \leq {(1-m) (N-2)}/{2}$. Hence
$$a(x,\tau')\, \Delta\tilde{B}^{\alpha}(x) - |\gamma|\, x \cdot \nabla\tilde{B}^{\alpha}(x) \le
\frac{m\, (|x|^2 + k_2)}{C^*} \Delta\tilde{B}^{\alpha}(x) -
|\gamma|\, x \cdot \nabla\tilde{B}^{\alpha}(x).$$ Again, by direct
computation (using also that $|\gamma| = \frac{1}{N-2-Nm}$ and
that $C^* = \frac{2m(N-2-Nm)}{1-m}$) we find
\begin{equation*}
\begin{split}
\frac{m\, (|x|^2 + k_2)}{C^*} &\Delta\tilde{B}^{\alpha}(x) - |\gamma|\, x \cdot \nabla\tilde{B}^{\alpha}(x) \\
&=- \frac{k_2\, (N-4-m(N-2))\, N}{2\, (N\, (1-m)-2)\, (k_2+|x|^2)}\, B^\alpha
= - \frac{\theta(m,n,k_2)}{ (k_2+|x|^2)^{\frac N2 - \frac 1{1-m}}} < 0
\end{split}
\end{equation*}
for $m < \frac{N-4}{N-2}$ and $\alpha = \frac{(N-2)(1-m)}{2} -
1$. From this the claim easily follows.
\\ \\
We will now  compare the terms $I_i$ in (\ref{eqn-limit}) in order
to get a strong contraction principle with the weight $B^{\alpha}$.

\noindent{\em Claim ii.} There is a uniform constant $C$
(independent of $R$) such that
\begin{equation}
\label{equation-L1-bound}
\int_{\mathrm{R}^N}|q(x,\tau)|\,  \eta_R  \tilde{B}^{\alpha}(x)  \, dx \le
\int_{\mathrm{R}^N} |\tilde{u}_0 -
\tilde{v}_0|\, \eta_R \tilde{B}^{\alpha}(x) \,dx + C.
\end{equation}
In particular, if $\tilde{u}_0 - \tilde{v}_0 \in
L^1(\tilde{B}^{\alpha},\mathrm{R}^N)$,  then $\tilde{u}(x,\tau) -
\tilde{v}(x,\tau) \in L^1(\tilde{B}^{\alpha},\mathrm{R}^N)$.

\noindent{\em Proof of Claim ii.}
The proof of Claim ii is similar to the proof of Lemma
\ref{lem-integrability}, once we know
\begin{eqnarray*}
I_i(R) &\le& C\int_{-\log T}^{\tau}\int_{R/2\le|x|\le R}|q|(x,\tau')\tilde{B}^{\alpha}(x)\, dx d\tau \\
&\le& C\int_{-\log T}^{\tau}\int_{R/2\le |x| \le R}\frac{dx}{|x|^{\frac{2(1+\alpha)}{1-m}+2}} \le C(\tau),
\qquad i\in\{1,2,3\},
\end{eqnarray*}
for our choice of $\alpha$, where $C$ is independent of $R$ (for
$I_2(R)$ see also the arguments in the proof of Claim iii). Estimate
\eqref{equation-L1-bound} now follows from the computation in
(\ref{equation-bounded-q}) and the claim i.

\noindent{\em Claim iii.} We have, $\lim_{R \to \infty} I_i(R)
=0$, $i=1,2,3$.
%there exists a uniform constant $C$, independent of $R$, so that
%$I_i \le C \qquad i = 1,2,3$.

\noindent{\em Proof of Claim iii.}
%By the comparison principle,
%$|q|(x,\tau) \le \tilde{B}_{k_2} - \tilde{B}_{k_1}$. That yields
%$$|q|\tilde{B}^{\alpha}(x,\tau) \le
%\frac{C}{|x|^{\frac{2(1+\alpha)}{1-m}+2}}, \qquad
%|q| \, |\nabla\tilde{B}^{\alpha}|(x,\tau) \le
%\frac{C}{|x|^{\frac{2(1+\alpha)}{1-m}+3}}$$ for $|x| >> 1$.
%Therefore,
%$$I_i \le C\int_{R/2 \le |x| \le
%R}\frac{dx}{|x|^{\frac{2(1+\alpha)}{1-m}+2}} \le \tilde{C}$$ for
%$\alpha\ge \frac{(N-2)(1-m)}{2} - 1$. Moreover, it is easy to see
%that $I_i \to 0$ as $R\to \infty$ in the case $\alpha >
%\frac{(N-2)(1-m)}{2} - 1$. This finishes the proof of the Claim.
%Choose $\alpha := \frac{(N-2)(1-m)}{2} - 1$ so that the above
%claims hold.
By the Claim ii  we have $|q|(\cdot,\tau') \in
L^1(\tilde{B}^{\alpha}, \R^N)$, for $\tau' \in [-\log
T,\tau]$. Hence, by choosing $R$ sufficiently big, the integral
$$\int_{-\log T}^{\tau}\int_{R/2 \le |x| \le
  R}|q|(x,\tau')\tilde{B}^{\alpha}(x)\,dx d\tau'$$ can be made
arbitrarily small. This readily implies that $\lim_{R \to \infty}
I_3(R)=0$. In addition, since $a(x,\tau) \leq C \, |x|^2$, this also
implies that $\lim_{R \to \infty} I_1(R) =0$. Finally, observing that
$|\nabla \tilde B^\alpha| \approx B^\alpha/R$ on $R/2 \le |x| \le R$,
we conclude that $\lim_{R \to \infty} I_2(R) =0$, finishing the proof
of the claim.

Combining the above two claims and \eqref{equation-bounded-q},  concludes the  proof of Lemma \ref{lem-contraction}.
\end{proof}

Since we have Lemma \ref{lem-contraction}, the following result holds
due to Osher and Ralston if we replace the usual $L^1$ norm by the
weighted $L^1(\tilde{B}^{\alpha}, \R^N)$ norm. This replacement will
leave the proof of the following Lemma unchanged.

\begin{lem}[Lemma 1 in \cite{OR} by Osher and Ralston]
\label{lem-or}
Let $R \ge R_0$. Suppose $\|\tilde{u}(\cdot,\tau_i) -
\tilde{v}_0\|_{L^1(\tilde{B}^\alpha, \R^N)} \to 0$ as $i\to\infty$, and let $\tilde{B}_k$
be any stationary solution of (\ref{eqn-rescaled}). If $\tilde{v}$
is a solution of (\ref{eqn-rescaled}) in $\R^N\times [0,\infty)$
with initial value $\tilde{v}(x,0) = \tilde{v}_0(x)$, then
$$\|\tilde{v}(\cdot,\tau) - \tilde{B}_k\|_{L^1(\tilde{B}^{\alpha}, \R^N)} =
\|\tilde{v}_0 - \tilde{B}_k\|_{L^1(\tilde{B}^{\alpha}, \R^N)},$$
for all $\tau > 0$ and all $k > 0$.
\end{lem}

\noindent{\em Remark on the proof of  Lemma \ref{lem-or}:} Define
$T(t)\tilde{u}_0 = \tilde{u}(t)$, where $\tilde{u}(t)$ is a
solution of (\ref{eqn-rescaled}) starting at $\tilde{u}_0$. The
proof of Lemma \ref{lem-or} uses only that $T(t)$ is a semi-group
on an $L^1(\tilde{B}^{\alpha}\eta_R,\R^N)$-closed subset of
$L^{\infty}$, satisfying the contraction principle
(\ref{eqn-contract}) and fixing $\tilde{B}_k$.

The following simple convergence result will be used in the proof
of Theorem \ref{thm-nonint}.

\begin{lem}
\label{lem-limit} Let $u_0$ satisfy (\ref{eqn-trapped}) for some
constants $k_1, k_2$.  Take any $\tau_i\to \infty$ and let
$\tilde{u}_i(\cdot, \tau) = \tilde{u}(\cdot, \tau_i + \tau)$.
Then,  passing to a subsequence, $\tilde{u}_i$ converges, as $i
\to\infty$,  uniformly on compact subsets of $\R^N \times (-\infty,\infty)$  to $\tilde{v}(x,\tau)$, an eternal
solution of (\ref{eqn-rescaled}), satisfying
(\ref{eqn-rescaled-trap}).
\end{lem}

\begin{proof}
Since $\tilde{u}$ satisfies (\ref{eqn-rescaled-trap}), equation
(\ref{eqn-rescaled}) is uniformly parabolic on $B_{2R}\times
[-(\log T)/2 - \tau_i,\infty)$, for any $R > 0$. By  standard
parabolic estimates, the sequence $\{\tilde{u}_i\}$ is equicontinuous on
compact subsets  of $\R^N \times (-\infty, \infty)$.  Hence, by  Arzela-Ascoli theorem and the
diagonalization argument any sequence $\{\tilde{u}_i\}$
will have a convergent subsequence, converging uniformly on compact
subsets to $\tilde{v}$, an eternal solution
of (\ref{eqn-rescaled}) satisfying (\ref{eqn-rescaled-trap}).
\end{proof}

Recall that we denote $\tilde{B}_{k_0}$ simply by $\tilde{B}$.

\begin{claim}
The sequence $\tilde{u}_i(x,\tau)$ converges to
$\tilde{v}(x,\tau)$ in $L^(\tilde{B}^{\alpha},\R^N)$ norm.
\end{claim}

\begin{proof}
By Lemma \ref{lem-contraction} and Lemma \ref{lem-contraction2} we have
$$\int_{\R^N}|\tilde{u}_i(x,\tau)-\tilde{B}(x)|\, dx \le C, \qquad \mbox{and}$$
$$\int_{\R^N}|\tilde{u}_i(x,\tau)-\tilde{B}(x)|\cdot |\tilde{B}^{\alpha}(x)\, dx \le C,$$
where $C$ is a constant independent of $i$ and $\tau$ and therefore,
\begin{eqnarray*}
& & \int_{\R^N}|\tilde{u}_i(x,\tau) -
\tilde{v}(x,\tau)|\tilde{B}^{\alpha}(x)\, dx \le \\
&\le& \int_{B^R}|\tilde{u}_i(x,\tau) -
\tilde{v}(x,\tau)|\tilde{B}^{\alpha}(x)\, dx + \\
&+& (\int_{|x| \ge R}|\tilde{u}_i(x,\tau) -
\tilde{B}(x)|\tilde{B}^{\alpha}(x)\, dx
+ \int_{|x| \ge R}|\tilde{v}(x,\tau) - \tilde{B}(x)|\tilde{B}^{\alpha}(x)\, dx) \\
&\le& \int_{B^R}|\tilde{u}_i(x,\tau) -
\tilde{v}(x,\tau)|\tilde{B}^{\alpha}(x)\, dx + \\
&+& \frac{C}{R^{2\alpha/(1-m)}}(\int_{|x| \ge
R}|\tilde{u}_i(x,\tau) - \tilde{B}(x)|\, dx
+ \int_{|x| \ge R}|\tilde{v}(x,\tau) - \tilde{B}(x)|\tilde{B}^{\alpha}(x)\, dx) \le \\
&\le& \int_{B^R}|\tilde{u}_i(x,\tau) -
\tilde{v}(x,\tau)|\tilde{B}^{\alpha}(x)\, dx +
\frac{C}{R^{2\alpha/(1-m)}}.
\end{eqnarray*}
If we let $i\to\infty$ in the previous estimate, since $\tilde{u}_i(x,\tau)\to \tilde{v}(x,\tau)$
uniformly on compact sets,  we get
$$\lim_{i\to\infty}\int_{\R^N}|\tilde{u}_i(x,\tau) - \tilde{v}(x,\tau)|\tilde{B}^{\alpha}(x)\,dx
\le \frac{C}{R^{2\alpha/(1-m)}},$$ which holds for every $R > 0$
and therefore $||\tilde{u}_i(\cdot,\tau) -
\tilde{v}(\cdot,\tau)||_{L^1(\tilde{B}^{\alpha}(x),\R^N)} \to 0$
as $i\to\infty$.
\end{proof}

\noindent{\em Proof of Theorem \ref{thm-nonint}.}
\begin{proof}
For a sequence $\tau_i \to \infty$, let $\tau_{i_k}$ a subsequence for
which $\tilde u(\cdot,\tau_i) \to \tilde v_0$, as $i_k \to\infty$,
uniformly on compact sets of $\R^N$, as shown in Lemma
\ref{lem-limit}. We will show that $\tilde v_0 = \tilde B$, as
stated in the Theorem.  By the previous claim and Lemma \ref{lem-or}
we have that
\begin{equation}
\label{equation-OR} ||\tilde{v}(\cdot,\tau) -
\tilde{B}(\cdot||_{L^1(\tilde{B}^{\alpha},\R^N)} =
||\tilde{v}(\cdot,\tau_0) -
\tilde{B}(\cdot)||_{L^1(\tilde{B}^{\alpha}, \R^N)},
\end{equation}
for all $\tau > -log T$. On the other hand, if $\max
|\tilde{v(x,0)} - \tilde{B}| > 0$ we have the strong
contraction principle (\ref{equation-min-max}), which contradicts
(\ref{equation-OR}). Therefore, $\tilde{v}(x,0) =
\tilde{B}(x)$.
\end{proof}

\end{section}

\begin{section}{Solutions that live longer}\label{sec-example}

In the previous sections  we established the vanishing profile of
solutions  $u$ of equation  \eqref{eqn-u} with  inital data $u_0$
trapped in between two Barenblatt solutions with the same
vanishing time $T$,  i.e. when \eqref{eqn-igc1} holds true. We
showed that  if   $u_0$ satisfies \eqref{eqn-trapped}, then $u$
vanishes at time $T$ and the rescaled solution $\tilde{u}(x,\tau)
= (T-t)^{-\beta}u(x\,  (T-t)^{\gamma}, t)$, with $\tau=1/(T-t)$,
converges,  as $\tau\to \infty$,  to a rescaled Barenblatt
solution $\tilde B$.

In this section we will show the condition \eqref{eqn-trapped} is
necessary for Theorems \ref{thm-integrable} and \ref{thm-nonint} to
hold true.  We will prove Theorem \ref{thm-example} which presents an
example of a class of initial conditions $u_0$ which satisfy the
growth condition
\begin{equation}\label{eqn-igc22}
u_0(x) = \left ( \frac{C^* \, T}{|x|^2}  \right )^{\frac 1{1-m}} (1 + o(1)), \qquad \mbox{as} \,\, |x| \to \infty
\end{equation}
with  $C^*$  given by \eqref{eqn-C*},  for which the solution $u$
of  \eqref{eqn-u}  with initial data $u_0$ vanishes at time $T^*
>T$. In  addition, we will show that the solution $u$ remains strictly positive for
$t < T^*$ and it satisfies, as $|x| \to \infty$, the growth
conditions $u(x,t) \approx   C(t) \, |x|^{- \frac 2{1-m}}$,  with $C(t) >0$
on $0 < t < T$ and $u(x,t) = O (|x|^{-\frac m{N-2}})$, on
$T < t <T^*.$ In particular, $u$ becomes integrable on  $t >T$ and its vanishing
behaviour is given by a class of self-similar solutions $\theta(x,t)$.

It is well known that the Barenblatt solutions   given by \eqref{eqn-baren}  are not the only self-similar solutions of
equation \eqref{eqn-u}. It was shown in \cite{PZ} that \eqref{eqn-u}
possesses   self-similar solutions of the form
$$\Theta(r,t) = (T-t)^\alpha \, f(\eta), \qquad \eta=\frac{r}{(T-t)^\theta}, \,\,
\alpha=\frac{1-2\theta}{1-m},$$
where the function $f$ is a
solution of an elliptic  non-linear eigenvalue problem with
eigenvalue $\theta$ satisfying the bound
$$ - \frac{m}{(1-m)\, N-2} < \theta < \frac 12,$$
$f'(0) = 0$ and $f(\eta) = O(\eta^{-(N-2)/m})$ as $\eta\to\infty$.
%For the critical value $m=(N-2)/(N+2)$,  $\theta = 0$.
The solution $f$ was shown to be unique apart from a scaling
due to the invariance of $f$ under the transformation
$f(\eta;\lambda)= \lambda^{-\frac 2{1-m}} \, f(\eta/\lambda;1)$ and satisfies
\begin{equation}\label{eqn-growth-eta2}
f(\eta;1) \approx \eta^{-\frac{N-2}{m}}, \qquad \mbox{as} \,\, \eta \to \infty.
\end{equation}
%with
%\begin{equation}\label{eqn-nu}
%\nu= \frac{N-2}{m}-\frac2{1-m}.
%\end{equation}
%The solution $\Theta$ satisfies
%$$\int_{\R^N} \Theta(\cdot,t) \, dx = O((T-t)^{\alpha+N\, \theta}) \to 0, \qquad \mbox{as} \,\,\, t \uparrow T.$$
It was shown in \cite{GP} that for any radially symmetric solution of \eqref{eqn-u}
with initial data satisfying the growth condition
$u_0(r) = O (r^{-\frac{N-2}{m}})$ as $r  \to \infty$, then
the vanishing behavior of $u$ is described by the self-similar
solutions $\Theta$, i.e. there exists $\lambda >0$ such that the
rescaled solution satisfies
$$(T-t)^{-\alpha} u(\eta \, (T-t)^\theta,t) \to f(\eta; \lambda)$$
uniformly in $\eta \geq 0$.

\smallskip
%Consider next the following initial data
%\begin{equation}
%\label{eqn-initial-asymp}
%u_0(x) = \left (\frac{C^*\, T}{|x|^2 +
%1} \right )^{\frac{1}{1-m}} \left (1 + \frac{A}{(|x|^2 + 1)^{\frac{\nu}2} }+ o(\frac{1}{(|x|^2 + 1)^{\frac{\nu}2}})\right )
%\end{equation}
%with $\nu$ as in \eqref{eqn-nu}.
%Let $u$ be the solution of  (\ref{eqn-u}) with initial data  $u_0$. Notice that
%we can not impose pointwise an upper bound on $u_0$ on $\R^N$ by any Barenblatt
%solution. We will show that $u$ has a different  vanishing behaviour than  this established in Theorems \ref{thm-integrable} and \ref{thm-nonint}, as stated next.

In  the proof of Theorem \ref{thm-example} we will use the following lemma,
which is also of independent interest.

\begin{lem}
\label{lem-pointwise}
Assume that $v$ is a solution of \eqref{eqn-u} on $\R^N
\times (0,T)$. Assume that $v_0 \leq f$ with $f \in L^1(\R^N)$ and
radially symmetric. Then, at time $t >0$ the solution $v$ satisfies
the bound $$v(x,t) \leq \frac{C}{|x|^{\frac{N-2}{m}}}, \qquad |x|
>1.$$
\end{lem}

\begin{proof}  We introduce  the potential function
 $$
 w(x,t) = \int_0^t  v^m (x,s)   \, ds
$$
which  satisfies the inequality
\begin{equation}\label{eqn-potential2}
\Delta w \geq - f, \qquad \mbox{on} \,\,\, \R^N.
\end{equation}
Let
$$Z(x) = \frac 1{N(N-2)\omega_N} \int_{\R^N} \frac 1{|x-y|^{N-2}} \, f(y) \, dy$$
denote  the Newtonian potential of $f$, so that we have
$$\Delta (w - Z) \ge 0.$$
Since $f$ is integrable and radially symmetric, there exists
a constant $C < \infty$ for which
\begin{equation}\label{eqn-estZ2}
Z(x) \leq \frac{C}{|x|^{N-2}}, \qquad \forall x \in \R^N.
\end{equation}
Indeed, this follows from the observation that for a radially symetric
$f$ the Newtonian potential of $f$  is also given by
$$Z(r) = \int_{r}^\infty \frac{1}{\rho^{N-1}}\int_{|y| \leq \rho} f(y) \, dy \, d\rho.$$
Similarly as in the proof of Lemma \ref{lem-integrability} we get
\begin{equation}
\label{eqn-decayw}
w(x)= \int_0^t v^m(x,s) \, ds  \le Z(x) \le \frac{C}{|x|^{N-2}}.
\end{equation}
We will now use this bound together with the Aronson-B\'enilan
inequality
\begin{equation}\label{eqn-ABI}
v_t \leq \frac 1{(1-m)\, t}  \, v
\end{equation}
to conclude the desired bound on $v$. Indeed, we first integrate
\eqref{eqn-ABI} in time to obtain the inequality
$$\frac{v(x,t_2)}{v(x,t_1)} \leq \left ( \frac {t_2}{t_1} \right )^{\frac 1{1-m}},
\qquad \mbox{if} \,\,  0 < t_1 < t_2.$$
Hence
$$w(x) \geq \int_{t/2}^t v^m(x,s) \, ds \geq v^m(x,t) \int_{t/2}^t
\left ( \frac {s}{t} \right )^{\frac m{1-m}} = c(t) \, v^m(x,t)$$
which combined with \eqref{eqn-decayw} implies the bound
$$v(x,t) \leq \frac{C(t)}{|x|^{\frac{N-2}{m}}}$$
as desired.
\end{proof}

We now proceed with the proof of the Theorem.

\noindent{\em Proof of Theorem \ref{thm-example}.} We begin with the following observation: If $w$ is a solution of \eqref{eqn-u} which vanishes at time
$T$, then for any $a >0$, the solution of \eqref{eqn-u} given by  $W(x,t) =
w(a\, x,a^2\, t)$ has vanishing time $T'=T/a^2$. Hence, we can make
$T'$ arbitrarily large by choosing $a$ sufficiently small.

Let $f \geq 0$, be any radially symmetric integrable function such that
$f(x) = o(|x|^{-\frac2{1-m}})$, as $|x| \to \infty$. In particular, we can take $f(x)$ to satisfy
$f(x) = O(|x|^{-\frac{N-2}m})$, $|x| \to \infty$.
Choose
$a$ sufficiently large so that the vanishing time $T'$ of the solution $w$
of \eqref{eqn-u} with initial data $w_0(x)=f(a  x)$ satisfies $T' >T$, with $T$
as in \eqref{eqn-igc22}.
Set $$u_0 = \left ( \frac{C^* \, T}{|x|^2 +1 }  \right )^{\frac 1{1-m}}  + w_0$$
and let
$B(x,t) = \left ( \frac{C^* \, (T-t) }{|x|^2 +1 }  \right )^{\frac 1{1-m}}$
the Barenblatt solution with   $B(x,0) = \left ( \frac{C^* \, T}{|x|^2 +1 }  \right )^{\frac 1{1-m}}.$
Denoting by $u$ the solution of \eqref{eqn-u} with initial data $u_0$ it is clear
application of the comparison principle that  $u\geq w$ so that the vanishing time $T^*$ of $u$ satisfies $T^* \geq T' >T$. This proves (i). Also, since $f \geq 0$, $u \geq B$, so that \eqref{eqn-igc3} is satisfied as well.

Since $u_0 - B(\cdot,0) = f \in L^1(\R^N)$, by Corollary
\ref{cor-contraction}
$$\|{u}(\cdot,t) - B(\cdot,t)\|_{L^1(\R^N)} \le \|f(\cdot)\|_{L^1(\R^N)}.$$
Since $B(\cdot,t) \equiv 0$ for $t \ge T$, this implies
\begin{equation}
\label{eqn-l1-finite}
\|{u}(\cdot,t)\|_{L^1(\R^N)} \le \|f\|_{L^1(\R^N)} \qquad
\mbox{for} \,\,\, t\in [T,T^*).
\end{equation}
Combining the estimate (\ref{eqn-l1-finite}) and Lemma \ref{lem-pointwise},
if we take  $v(x,0) = {u}(x,T)$, yields to \eqref{eqn-igc4}.

The statement (iii) of our
Theorem now immediately follows by the result of Galaktionov and Peletier in \cite{GP}. \qed

The proof of Corollary \ref{cor-yamabe} now easily follows.

\noindent{\em Proof of Corollary \ref{cor-yamabe}.}
It is known that if $m = \frac{N-2}{N+2}$ then $\alpha = \frac{N+2}{4}$ and
$\theta = 0$ in the definition of self-similar solutions $\Theta(x,t)$. In this case
function $f$ is given explicitly by 
$$f(\eta,\lambda) = \left (\frac{K_N \, \lambda}{\lambda^2 + \eta^2} \right )^{\frac{N+2}{2}}.$$
The proof  of Theorem \ref{thm-example} implies  that for  $t > T$, we have the bound 
$$u(x,t) \le \frac{C(t)}{|x|^{N+2}}.$$
By the result of Del Pino and Saez in \cite{DS}, if the vanishing time of $u$ is $T^*$ with  $T^* > T$, then 
there is $\lambda > 0$ so that
$$(T^* - t)^{-\frac{1}{1-m}}\, u(x,t) \to f(|x|,\lambda) \qquad \mbox{as}\,\,\, t\to T^*.$$
\qed

We will end the paper  with the following remark regarding the examples of solutions $u$  constructed in
Theorem \ref{thm-example}.  We know that these solutions $u$   in  live up to $T^* > T$,
and satisfy $u(x,t) \ge B(x,t)$ for $t < T$ and $u(x,t) \le
{C(t)}\, {|x|^{-\frac{N-2}m}}$ for $t > T$.

\begin{remark}
In the examples of solutions  $u$ constructed in Theorem \ref{thm-example}, the  Barenblatt $B(x,t)$ dies off exactly at time $T$,  in the sense that  there is
no sequence $|x_i|\to \infty$  and $A > 0$ so that
\begin{equation}
\label{equation-big} u(x_i,T) \ge
\frac{A}{|x_i|^{\frac{2}{1-m}}}.
\end{equation}
\end{remark}

\begin{proof}
Lets argue by contradiction. Assume there is a sequence
$|x_i|\to\infty$ so that (\ref{equation-big}) holds for some $A
> 0$. Take $t > T/2$ so that $2^m(C^*(T-t))^{\frac{m}{1-m}} < A/2$
and let $\delta = (T-t)/2$. Notice that by our choice of $f$, $0 \leq
u(x,0) - B(x,0) \in L^1(\R^N)$. By the comparison principle, we have   $u(x,t) \ge B(x,t)$ for $t > 0$,  and by Lemma
\ref{lem-integrability}, $u(x,t) - B(x,t) \in L^1(\R^N)$. Let
$$w(x,t) = \int_t^{t+\delta}(u^m - B^m)\,ds.$$
As before, $\Delta w \ge -(u - B)$ and since $||(u -
B)(\cdot, t)||_{L^(\R^N)} \le ||f||_{L^1(\R^N)}$, using  Newtonian
potentials we obtain that 
$$w(x,t) \le \frac{C}{|x|^{N-2}}$$
for a uniform constant $C$. This yields the existence of $s\in
(t,t+\delta)$ such that
\begin{equation}
\begin{split}
\label{equation-vanishing}
u^m(x,s) &\le B^m(x,s) + \frac{C}{\delta \, |x|^{N-2}} \\
&\le
\left ( \frac{C^* \, (T-s)}{|x|^2}\right )^{\frac{m}{1-m}}
+ \frac{C}{\delta\cdot |x|^{N-2}} . 
\end{split}
\end{equation}
For each $x_i$ choose $s_i\in [t, t+\delta]$ so that
(\ref{equation-vanishing}) is satisfied. The Aronson-B\'enilan
inequality gives
$$u(x_i,T) \le \frac{T}{s_i}\, u(x_i,s_i) \le 2\, u(x_i,s_i).$$
Combining (\ref{equation-big}) and (\ref{equation-vanishing})
yields
$$\frac{A}{|x_i|^{\frac{2m}{1-m}}} \le u^m(x_i,T) \le \frac{2^m(C^*(T-s_i))^
{\frac{m}{1-m}}}{|x_i|^{\frac{2m}{1-m}}} + \frac{C}{\delta \, 
|x_i|^{N-2}}$$ and therefore
$$\frac{A}{2\, |x_i|^{\frac{2m}{1-m}}} \le \frac{C}{\delta\,  |x_i|^{N-2}}$$
which can not be fulfilled for $|x_i| >> 1$. This finishes the proof of our claim. 
\end{proof}

%Almost the same proof as that of the Theorem \ref{thm-example} yields
%the following result.

%\begin{cor}
%Let $u_0(x) = (\frac{T}{C^*(|x|^2 + k)})^{\frac{1}{1-m}}(1 + o(1))$, where
%$o(1) = f(x) \ge 0$ and $f(x) \le h(|x|)\in L^1(\R^N)$ (where $h(r)$ is a radial
%function). Then,
%\begin{enumerate}
%\item[(i)]
%either $u(x,t)$ vanishes at $T$, or
%\item[(ii)]
%$u(x,t)$ vanishes at $T^* > T$ and for $t > T$ the bahaviour changes from
%$u(x,t) = O(\frac{1}{(|x|^2 + 1)}^{1/(1-m)})$,
%as $|x| \to\infty$, to that of the finite area. In particular,
%$$(T^* - t)^{-\frac{1}{1-m}}u(x,t) = \{\frac{k_N}{1 + |x|^2}\}^{\frac{N+2}{2}}
%+ \theta(x,t), \qquad \mbox{for} \,\,\, t > T,$$
%and $\sup_{x\in \R^N}(1 + |x|^{N+2})|\theta(x,t)| \to 0$ as $t\to T^*$
%\end{enumerate}
%\end{cor}

\end{section}

\section{Appendix}

In this appendix we improve Theorem \ref{thm-nonint}. The goal is to remove
the assumption  $u_0(x) \geq  B_{k_1}(x,0)$, where $B_{k_1}(x,t)=({C^*\, (T-t)}/(k_1 + |x|^2))^{\frac{1}{1-m}} $ is a Barenblatt solution with the same vanishing time as $u$. 
This was assumed in \eqref{eqn-trapped}. The bound from above is necessary, as was proven in the previous section. Denoting by  $B(x,t)=({C^*\, (T-t)}/(k + |x|^2))^{\frac{1}{1-m}} $, for some $k >0$, we will  show   the following result.

\begin{thm}
\label{thm-profile-bigger-6}
If $0 < m \le \frac{N-4}{N-2}$ for $N > 4$ and
\begin{equation}
\label{equation-cond2}
|u_0(x) - B(x,0)| \le f(|x|) \in L^1(\R^N) \quad \mbox{and} \quad  u_0(x) \le B(x,0)
\end{equation}
for a positive, radial function $g$, then $u$ vanishes at the same
time as $B(x,t)$ and the rescaled solution $\tilde{u}(x,\tau)$
converges, as $\tau\to\infty$, uniformly $\R^N$, to the rescaled
Barenblatt $\tilde{B}$.
\end{thm}

%Theorem \ref{thm-profile-bigger-6} removes the bound from below by a
%Barenblatt solution that we have assumed in Theorem \ref{thm-nonint}.
%We will prove Theorem \ref{thm-profile-bigger-6}
%in the sequa of lemmas and claims. 

%\begin{proof}[Proof of Theorem \ref{thm-profile-bigger-6}] 
To simplify the notation we will assume that $B(x,t)=({C^*\, (T-t)}/(1 + |x|^2))^{\frac{1}{1-m}} $. The proof of the theorem will be based
on a sequence of  observations. 

We begin by 
noticing  that condition (\ref{equation-cond2})  implies the
$L^1$-contraction principle
\begin{equation}
\label{eqn-contract}
\int_{\R^N}|u(\cdot,t) - B(\cdot,t)|\, dx \le \int_{\R^N}|u_0 - B(\cdot,0)|\, dx,
\qquad \forall  t\in [0,T).
\end{equation}
The proof is the same as the proof of Corollary 
\ref{cor-contraction}, which goes through under the  weaker assumption that only $v_0 \geq B$ (and not necessarily $u_0 \geq B$).  Furthermore, this implies
that $u$ and $B$ have the same vanishing time.

%\begin{cor}
%If $u(x,t)$ and $B_k(x,t)$ are as above, then the solutions $u(x,t)$ and $B_k(x,t)$
%have the same vanishing time.
%\end{cor}

%\begin{proof}
%Let $T^*$ and $T$ be the vanishing times of $u$ and $B_k$ respectively. By the comparison
%principle, $T^* \le T$. If  $T^* < T$,  (\ref{eqn-contract}) would imply
%$B(x,t) \in L^1(\R^N)$ for $t \in [T,T^*)$ which is not true. 
%\end{proof}
\smallskip

In order to be able to take the limit of the rescaled solution, we
need to establish the necessary   a'priori estimates. It turns out that it is possible
to do so, just by using the fact that the difference $u(\cdot,0) -
B(\cdot,0) \in L^1(\R^N)$ and that $B$ and $u$ have the same vanishing
time.  Introducing as in the previous sections the rescaling 
\begin{equation}
\label{equation-rescaling}
\tilde{u}(x,\tau) = (T-t)^{-\beta}u(x\, (T-t)^{\gamma},\tau), \qquad \tau = -\ln(T-t)
\end{equation}
(with $\beta, \gamma$ given by \eqref{eqn-parameters}),  which satisfies the rescaled equation 
\begin{equation}
\label{equation-rescaled}
(\tilde{u})_{\tau} = \Delta\tilde{u}^m + |\gamma|\, \dv(x\cdot\tilde{u}), \qquad
\mbox{on} \,\, \R^N \times (-\log T,\infty)
\end{equation}
we have:

\begin{prop}
\label{prop-bounds}
If $|u_0(x) - B(x,0)|\le f(|x|) \in L^1(\R^N)$, then there are
positive constants $C_1, C_2, r_0, \tau_0$ such  that
\begin{equation}
\label{equation-good-bounds}
\frac{C_1}{(r^2 + 1)^{1/(1-m)}} \le \tilde{u}(r,\tau) \le \frac{C_2}{(r^2 + 1)^{1/(1-m)}} \qquad \mbox{for} \,\, r\ge r_0 
\end{equation}
for all $\tau \in [\tau_0,\infty)$.
\end{prop}

\begin{proof}
We will first prove \eqref{equation-good-bounds} under the assumption that $u_0$
is radially symmetric. At the end of the proof we will remove this assumption.

As we have observed in the proof
of Lemma \ref{lem-integrability},  absolute value $|u-B|$ satisfies the differential inequality
$$\frac{d}{dt}\, |u - B| \le \Delta|u^m - B^m|.$$
Hence, for any fixed numbers  $T/2 < s < t < T$, the  function  
$$w(x) = \int_s^t |u^m - B^m|(x,l) \, dl$$  satisfies 
$$\Delta w \ge - |u - B|(s).$$
Let $Z$ be such that $\Delta Z = -|u-B|(s)$ and therefore
$$\Delta(w - Z) \ge 0$$
which, by the mean value property, as in the proof of  Lemma \ref{lem-integrability}, yields to the inequality
$$w(x) \le Z(x,s),\qquad r=|x|.$$
Since, $u-B$ is radially symmetric, the function potential $Z$ is given by
$$Z(x,s) = c_n\,  \int_r^{\infty}\frac{1}{\rho^{N-1}}\int_{|y|\le\rho}|u-B|(s)\,dy\,d\rho, \qquad r=|x|$$
for an appropriate constant $c_n$. 
We conclude that
$$w(x) \le \frac{||(u-B)(s)||_{L^1}}{r^{N-2}}, \qquad r \geq 1$$
which combined with \eqref{eqn-contract} and our assumption \eqref{equation-cond2}
implies the bounds
$$
\int_s^t B^m(x,l)\, dl  - C \, \frac{\|f\|_{L^1}}{r^{N-2}} \le
\int_s^t u^m(x,l)\, dl  \le \int_s^t B^m(x,l) \, dl+ C \, \frac{\|f\|_{L^1}}{r^{N-2}}.
$$
Next, fix $t \in [3T/4, T)$ and choose $s \in [T/2,T)$ such that $T-t=t-s$, so that 
$$2^{-\frac m{1-m}}\, (t-s) \, B^m(x,s) \leq  \int_s^t B^m (x,l) \, dl \leq 2^{\frac m{1-m}}\, (t-s) \, B^m(x,t).$$
Combining the two last inequalities gives
\begin{equation*}
2^{-\frac m{1-m}}\,B^m(x,s)  - C \, \frac{\|f\|_{L^1}}{r^{N-2}} \leq \frac 1{t-s} \,\int_s^t u^m(x,l)\, dl 
\leq  2^{\frac m{1-m}}\,B^m(x,t)  + C \, \frac{\|f\|_{L^1}}{r^{N-2}}.
\end{equation*}
>From the above we conclude a pointwise bound from above and below for
$u^m(x,t)$ with the aid of the  Aronson-Benilan inequality, $(\log u)_t \le 1/ ((1-m)\, t)$,
which after integration implies the bound 
$$(\frac st)^{\frac m{1-m}} \, u^m(x,t) \leq  \frac 1{t-s}  \int_s^t  u^m(x,l)\, dl \leq (\frac ts)^{\frac m{1-m}} \, u^m(x,s).$$
Due to our choice for $s$ (for a given $t$) we have $t/s \leq 2$. Hence, combining the
last two inequalities gives
\begin{equation}
\label{equation-less}
u^m(r,t) \le  C_1\, B^m(r,t) +
\frac{C}{(T-t)\, r^{N-2}}
\end{equation}
which holds for all $t\in [3T/4,T)$, since $t$ is arbitrary, and also
\begin{equation}
\label{equation-bigger}
u^m(r,s) \ge C_2 \, B^m(r,s) -
\frac{{C}}{(T-s)\, r^{N-2}}
\end{equation}
which also holds for all  $s\in [3T/4,T)$  since $t$ is arbitrary and $T-t=t-s$. 
Rescaling inequalities \eqref{equation-less} and \eqref{equation-bigger}, 
and using \eqref{eqn-parameters} we conclude 
%$$\frac{C_1}{(r^2 + 1)^{m/(1-m)}} - \frac{\tilde{C}_1\cdot(T-t)^{-m\beta}}{r^{N-2}\cdot (T-t)^{(N-2)\gamma}+1}
%\le \tilde{u}^m(r,t) \le
%\frac{C_2}{(r^2 + 1)^{m/(1-m)}} + \frac{\tilde{C}_2\cdot(T-t)^{-m\beta}}{r^{N-2}\cdot (T-t)^{(N-2)\gamma + 1}},$$
$$\frac{C_1}{(r^2 + 1)^{m/(1-m)}} - \frac{C}{r^{N-2}} \le u^m(r,\tau)
\le \frac{C_2}{(r^2 + 1)^{m/(1-m)}} + \frac{C}{r^{N-2}}$$
for some uniform constants $C_1, C_2, C$  and $\tau \geq \tau_0:=- \log (T/4)$. 
This readily implies \eqref{equation-good-bounds}
for $r \geq r_0$ (independent of $\tau$)  if $N-2 > \frac{2m}{1-m}$.
The last is equivalent to  $m < \frac{N-2}{N}$ and is implied by our assumption
$m \leq (N-4)/(N-2)$. 

In the case where $u(x,0)$ is nonradial and $B(x,0) - u(x,0)$ is bounded
from above by a radial function in  $ L^1(\R^N)$, define
$$\underline{u}_0(r) := \inf_{|x| =r}u(x,0) \quad \mbox{and} \quad  \overline{u}_0(r) := \sup_{|x|=r}u(x,0).$$
By our assumption, $B(x,0) - f(r) \le u(x,0) \le B(x,0) + f(r)$, where $r=|x|$ 
and $f\in L^1(\R^N)$. Since $B(x,0)$ is a radial function itself, we have 
$B(x,0) - \underline{u}_0(r) \in L^1(\R^N)$ and $B(x,0) - \overline{u}_0(r) \in L^1(\R^N)$.
Let $\underline{u}(x,t)$ and $\overline{u}(x,t)$ be the solutions to (\ref{eqn-u}) with initial data $\underline{u}_0(r)$ and $\overline{u}_0(r)$, respectively.
Then,  the radial result implies the bounds 
\begin{equation*}
\frac{C_1}{(r^2 + 1)^{m/(1-m)}} - \frac{C}{r^{N-2}} \le \underline{u}^m(r,\tau)
\le \frac{C_2}{(r^2 + 1)^{m/(1-m)}} + \frac{C}{r^{N-2}},
\end{equation*}
and
\begin{equation*}
\frac{C_1}{(r^2 + 1)^{m/(1-m)}} - \frac{C}{r^{N-2}} \le \overline{u}^m(r,\tau)
\le \frac{C_2}{(r^2 + 1)^{m/(1-m)}} + \frac{C}{r^{N-2}}.
\end{equation*}
Since, by  the comparison principle, $\underline{u}(r,\tau) \le u(x,\tau) \le \bar{u}(r,\tau)$, the above inequalities imply  \eqref{equation-good-bounds} in the nonradial case.
\end{proof}

\begin{proof}[Proof of Theorem \ref{thm-profile-bigger-6}]  
The previous Proposition yields the existence of $r_0 > 0$ and $ \tau_0 < \infty$ so that 
$$\frac{C_1}{(r^2+1)^{1/(1-m)}} \le \tilde{u}(r,\tau) \le \frac{C_2}{(r^2+1)^{1/(1-m)}}$$
for $\tau \in [\tau_0,\infty)$ and $r\ge r_0$. Let $Q_{r_0} = B_{r_0}\times [\tau_0,\infty)$.
Hence,  there exists a constant $c_0=c(r_0,\tau_0) >0$ such that 
$$\tilde{u}(r,\tau) \ge \frac{C_0}{(r_0^2 + 1)^{1/(1-m)}}, \qquad \mbox{on}\,\,  \partial_p Q_{r_0}.$$
By the maximum principle
$$\inf_{Q_{r_0}}\tilde{u} \ge \frac{C_0}{(r_0^2 + 1)^{1/(1-m)}}$$
which, combined with the lower bound in \eqref{equation-good-bounds}  implies that
$$\tilde{u}(x,\tau) \ge \frac{C_1}{(r^2+1)^{1/(1-m)}}, \qquad \mbox{on}\,\,\R^N\times [\tau_0,\infty) $$
 for a constant $C_1$ that depends on $r_0$. 

By our assumption we have
$\tilde{u}(x,\tau) \le \tilde{B}(x)$ on $\R^N\times [-\log T,\infty)$.  Hence, there are uniform 
constants $C_1$ and $C_2$ such  that
\begin{equation}
\label{equation-bounds-all-space}
\frac{C_1}{(r^2+1)^{1/(1-m)}} \le \tilde{u}(x,\tau) \le \frac{C_2}{(r^2+1)^{1/(1-m)}},
\qquad \mbox{on}\,\,\R^N\times [\tau_0,\infty).
\end{equation}
We conclude that the difference $\tilde u - \tilde B$ satisfies the equation
$$(\tilde{u}-\tilde{B})_{\tau} = \Delta(\tilde{a}(\tilde{u}-\tilde{B})) + |\gamma|\, \dv(x\cdot(\tilde{u}-\tilde{B})),$$
with $\tilde{a}(x,\tau) = \int_0^1\frac{d\theta}{(\theta\tilde{u} + (1-\theta)\tilde{B})^{1-m}}$ satisfying the bounds 
\begin{equation}
\label{equation-ellipticity}
\tilde{C}_1\, (r^2 + 1) \le \tilde{a}(x,\tau) \le \tilde{C}_2\, (r^2 + 1)
\end{equation}
on $R^N\times [\tau_0,\infty)$.
The rest of the proof of Theorem \ref{thm-profile-bigger-6} is the same as that of
Theorem \ref{thm-nonint}.
\end{proof}

\end{document}